\documentclass[11pt,twoside]{article}
\usepackage{mathrsfs}  
\usepackage{amsmath}   
\usepackage{amssymb}   
\usepackage{amsthm}    
\usepackage{esint}     
\usepackage{mathtools} 
\usepackage{fancyhdr}  
\usepackage{latexsym}  
\usepackage{bbding}    
\usepackage{wasysym}   
\usepackage{cite}      
\usepackage{authblk}
\usepackage{color}
\usepackage{multicol,graphics}
\newtheorem{theorem}{Theorem}[section]
\newtheorem{lemma}{Lemma}[section]
\newtheorem{proposition}{Proposition}[section]
\newtheorem{definition}{Definition}[section]
\newtheorem{remark}{Remark}[section]

\newtheorem{problem}{Problem}
\newenvironment{Proof}[1][Proof]{\noindent\textbf{#1.} }{\hfill $\Box$}

\makeatletter
\allowdisplaybreaks  
\begin{document}
\title{{On the Riemann Boundary Value Problem for Poly- and Meta-hyperanalytic Function Spaces over $d$-summable Curves}}
\author[a]{Yan Dai\thanks{\text{E-mail}: dygsmail0@csu.edu.cn}}
\author[b]{Juan Bory-Reyes\thanks{\text{E-mail}: juanboryreyes@yahoo.com}}
\author[a]{Fuli He\thanks{\text{Corresponding author. E-mail}: hefuli999@163.com}}
\affil[a]{School of Mathematics and Statistics, Central South University, Changsha 410083, China}
\affil[b]{ESIME-Zacatenco, Instituto Polit\'{e}cnico Nacional, M\'{e}xico, CDMX, 07738 M\'{e}xico}
\date{}
\maketitle
\begin{abstract}
Hyperanalytic functions, in the sense established by the mathematician Avron Douglis, are Douglis algebra-valued functions defined via a hypercomplex structure rather than the standard Cauchy-Riemann equations characteristic of traditional complex analysis. The classes of poly-hyperanalytic and meta-hyperanalytic functions represent advanced generalizations of Douglis's analysis. They are employed in the study of partial differential equations and elasticity, extending the concept of the classical holomorphic function through higher-order iterations and non-homogeneous terms.
The aim of this work is to find solvability conditions for a fundamental Riemann-type boundary value problem for spaces of poly-hyperanalytic and meta-hyperanalytic functions defined on an open, bounded, simply connected subset of the complex plane, where the boundary need only be a closed $d$-summable curve. In fractal geometry, $d$-summability is a geometric property—used to define the boundaries of fractal domains, enabling advanced mathematical integration and calculus on complex structures—defined by Jenny Harrison and Alec Norton. 
\end{abstract}

\noindent{\qquad\!\bf Keywords.} Douglis analysis, Riemann boundary value problem, non-rectifiable curve, poly- and meta-hyperanalytic functions.

\noindent{\qquad\!\bf AMS Subject Classifications.}  28A80, 30E20, 30E25, 30G35.

\section{Introduction}

The classical Riemann boundary value problem is a fundamental problem in the theory of analytic functions and singular integral equations, with significant applications in complex analysis, mathematical physics, and engineering. Specifically, it finds applications in plane elasticity (see e.g., \cite{Lu1}). It is also related to singular integral equations and can, in a certain sense, determine whether a singular integral equation is solvable (see e.g., \cite{MR1215485}). It was first formulated and studied by the German mathematician Bernhard Riemann \cite{Riemann} in the mid-19th century as part of his investigations into complex functions and differential equations. It is a type of boundary value problem for analytic functions, which aims to find sectionally analytic functions that satisfy certain internal and external boundary value relations on a given smooth closed curve. Using Cauchy-type integrals and the Plemelj formula, the standard works on Riemann boundary value problem are \cite{MR198152, MR1279172}.

Generalizations of the Riemann boundary value problem have been investigated together with new theoretical results not only for domains with non-smooth boundaries, but also for linear
and nonlinear elliptic systems in the plane, see \cite{GilbertBuchanan1983, Begehr94, Wendland}.

Let us indicate briefly some of the recent developments on the generalized Riemann boundary value problem. In \cite{Dyn79} a deeper discussion of the classical Riemann boundary value problem on rectifiable Jordan closed curves was presented. In 1983, the solvability of the Riemann boundary value problem on non-rectifiable closed curves was studied by Kats in \cite{Kat83}, and later by Liu in \cite{Liu}. In 2011, Abreu Blaya, Bory Reyes, and Kats generalized Cauchy-type integrals on non-rectifiable curves in \cite{MR2786194} and derived a distributional approach to the Riemann boundary value problem on non-rectifiable curves. In 2014, the authors in \cite{MR3146343} further improved the box dimension of curves to $h$-summability, refined the jump function to the generalized H\"{o}lder class, and studied and solved the Riemann boundary value problem on $h$-summable curves. As a generalization, this result enables the derivation of the solvability of the Riemann boundary value problem on \(d\)-summable curves \cite{greenthm, Geometintegrat}. While \cite{DBH} studies the Riemann boundary value problem on non-rectifiable curves for null-solutions of iterated and meta class of certain type of Beltrami equation. Moreover, the Riemann boundary value problem can also be discussed in high-dimensional spaces. In \cite{MR2302870}, Abreu Blaya, Bory Reyes, and Pe\~na Pe\~na generalized Kats's results to \(n\)-dimensional Euclidean space and solved the jump problem on general non-smooth surfaces.

The class of solutions of the classical Riemann boundary problem is the class of sectionally holomorphic functions, which can be generalized to the class of poly-analytic functions, the class of meta-analytic functions, the class of generalized analytic functions, and null-solutions of various elliptic systems. Poly-analytic functions arose from problems in planar elasticity, with Kolosov conducting the first investigations into them in 1908. Balk's remarkable monograph \cite{MR1184141} offers a comprehensive survey of poly-analytic and meta-analytic functions. In \cite{MR1979529} and \cite{MR2155424}, Du and Wang conducted a detailed discussion on the classical Riemann boundary value problem for poly-analytic and meta-analytic functions on simple smooth curves and the real axis. For the boundary value problems of these two types of functions, one can also refer to \cite{KHW, MR4165525, MR4596297}. A systematic and comprehensive study on boundary value problems with generalized function classes for the solution set can be found in \cite{MR150320}.

Hyperanalytic functions, introduced by Douglis in \cite{douglis1953}, are solutions of a special elliptic system of first-order partial differential equations in two independent variables. This system can be represented by a single hypercomplex equation involving the so-called Douglis operator. The theory of hyperanalytic functions, also referred to as Douglis analysis, generalizes the theory of analytic functions to the hypercomplex setting. For a thorough treatment of this theory, the reader is directed to \cite{MR3270573, MR3695391, GilbertBuchanan1983, Soldatov, Gilbert, GilbertHile, LiuChunGen, S}. Hyperanalytic functions have found numerous applications in plate and shell theories (see e.g., \cite{GilbertHackl, BuchananGilbert}).

The Riemann boundary value problem for hyperanalytic functions on rectifiable curves was studied in \cite{MR2005hyperRBVP, GilbertZeng1992, Zeng1991a, He}. In \cite{MR3800978, Zeng1990}, a direct generalization of the hyperanalytic Riemann boundary value problem for the non-rectifiable framework was given. In 2010, Abreu Blaya, Bory Reyes, and Vilaire considered the Riemann boundary value problem for hyperanalytic functions on $d$-summable curves (see \cite{Aa32}). Their work not only derived the solutions to the Riemann boundary value problem for hyperanalytic functions on general curves, but also provided a standard method for the Riemann boundary value problem for hyperanalytic functions on non-rectifiable curves, and offered strong inspiration for this paper. Furthermore, this paper generalizes the class of functions from hyperanalytic to poly-hyperanalytic and meta-hyperanalytic, thereby supplementing the results of \cite{Aa32}.

There has been considerable effort in the literature to formulate boundary problems of mathematical physics for domains with highly irregular boundaries, where it turns out that irregular media of fractal type are a good model. In this paper, we consider the class of \(d\)-summable curves---a mild modification of \(h\)-summable curves, which is a concept in fractal geometry describing curves that possess a specific type of geometric regularity, defined by means of the so-called gauge (dimension) functions; different functions of the diameter (such as \(x^d\)) give rise to \(d\)-summability.

The aim of this work is to derive solvability conditions and explicit solutions for Riemann boundary value problems associated with poly-hyperanalytic and meta-hyperanalytic functions with an identical factor on closed $d$-summable Jordan curves. 

The paper is structured as follows: After this brief introduction, Section 2 contains the relevant background information for stating and proving the main results. In particular, we introduce the Douglis algebra, hyperanalytic functions, and elaborate on the method by which Abreu Blaya, Bory Reyes, and Vilaire address the Riemann boundary value problem on \(d\)-summable curves in \cite{Aa32}. In Section 3, we present the Riemann boundary value problem for poly-hyperanalytic functions with identical factors on $d$-summable Jordan closed curves, and derive its solvability conditions as well as explicit solutions. In Section 4, we present the Riemann boundary value problem for meta-hyperanalytic functions with identical factors on $d$-summable Jordan closed curves, along with its solvability conditions and explicit solutions.

\section{Preliminaries}

In this section, we first introduce some basic concepts and notations of curves, Douglis algebras, hyperanalytic functions, and function classes required for Riemann boundary value problems.

\subsection{Curves and function classes}
A curve $\Gamma$ is the image of a continuous map $\gamma$ from the closed bounded interval $[0,1]$ (or $[a,b]$) to $\mathbb{C}$.
\begin{definition}
We call \( \Gamma \) rectifiable if \( \gamma \) is a function of bounded variation (see e.g., \cite{MR503901}). The length of the curve \( \Gamma \) is defined to be the total variation of \( \gamma \).
\end{definition}
\begin{definition}
The curve \( \Gamma \) is called a Jordan closed curve if \( \gamma(0)=\gamma(1) \) and \( \gamma \) is injective on \( (0,1) \). We say that \(\Gamma\) is a simple smooth closed curve if \( \Gamma \) is a Jordan closed curve, \( \gamma \) has continuous derivative on \( [0,1] \), $|\gamma'|$ does not vanish on $[0,1]$, and $\gamma'_{+}(0)=\gamma'_{-}(1)$.
\end{definition}
\begin{remark}
A Jordan closed curve may be rectifiable, such as a simple smooth closed curve, and may also be non-rectifiable, such as a fractal curve (see e.g., \cite{MR1102677}).
\end{remark}
The interior of the Jordan closed curve \( \Gamma \) is denoted by \( D^+ \), which is a bounded connected open set. And the exterior of \( \Gamma \) is denoted by \( D^- \), which is an unbounded connected open set. We have that \( \mathbb{C}\setminus\Gamma=D^+\cup D^-\). In the following text, unless explicitly stated otherwise, both \( D^+ \) and \( D^- \) are defined in this manner. We usually call \(\Gamma\) the jump curve. For complex functions, we can define their boundary values on the jump curve \( \Gamma \) (see e.g., \cite{MR1215485,MR198152,MR1279172}).
\begin{definition}
Let \(\Gamma\) be a Jordan closed curve. Suppose that \(F : \mathbb{C}\setminus\Gamma\to\mathbb{C}\) is a complex valued function and denote \(\lim\limits_{z\to t,z\in D^{+}}F(z)\) (resp. \(\lim\limits_{z\to t,z\in D^{-}}F(z)\)) by \(F^+(t)\) (resp. \(F^-(t)\)), if \(\lim\limits_{z\to t,z\in D^{+}}F(z)\) (resp. \(\lim\limits_{z\to t,z\in D^{-}}F(z)\)) exists for the fixed point \(t\in\Gamma\). \( F^\pm(t) \) are called boundary value functions.
\end{definition}
We often need to discuss the properties of complex functions at infinity, and for this purpose, we introduce the following definition (see e.g., \cite{MR2155424}).
\begin{definition}
Let \( F \) be defined on a neighbourhood of infinity, \( F \) is said to be of order \( m \) at infinity, denoted by \( \mathrm{Ord}(F,\infty)=m\), if
\[ \limsup\limits_{z\to\infty} |z^{-m}F(z)|=a>0.\]
\end{definition}
\begin{remark}
\( \mathrm{Ord}(F,\infty)=+\infty \), if for each \( m\in\mathbb{Z} \), \( \limsup\limits_{z\to\infty} |z^{-m}F(z)|=\infty \). And \( \mathrm{Ord}(F,\infty)=-\infty \) means that \( F(z)\equiv0 \). If \( F \) is a holomorphic function defined on a neighbourhood of the infinity, then \( F \) can be expanded in a Laurent series and \( \limsup\limits_{z\to\infty} |z^{-m}F(z)|=\lim\limits_{z\to\infty} |z^{-m}F(z)|= |a_m| >0 \) from \( \mathrm{Ord}(F,\infty)= m \) for \( m\in \mathbb{Z} \), which means that Laurent coefficients \( a_k = 0 \) for \( k>m \) and \( a_m\neq0 \).
\end{remark}

Next, we define an important class of functions.
\begin{definition}
Let \(A\) be a subset of \( \mathbb{C} \) and \( \nu\in(0,1] \). The class of H\"{o}lder functions on the set \( A \) is defined as
\[H^\nu(A)=\{f(t) : f \text{ is bounded and } p(f,\nu):=\sup\limits_{t,t'\in A,t\neq t'}\dfrac{|f(t)-f(t')|}{|t-t'|^\nu}<+\infty\}.\]
\end{definition}
\begin{remark}\label{properties_of_Holder_class}
It is easy to see that \(H^\nu(A)\) is a linear space and \( f \) is continuous for \( f\in H^\nu(A) \). When \( A \) is a compact set, then
\[H^\nu(A)=\{f(t) : p(f,\nu):=\sup\limits_{t,t'\in A,t\neq t'}\dfrac{|f(t)-f(t')|}{|t-t'|^\nu}<+\infty\},\]
and we have that \(f_1f_2\in H^\nu(A) \) for \( f_1,f_2 \in H^\nu(A) \) and \( H^{\nu_1}(A) \subseteq H^{\nu_2}(A) \) for \( 0<\nu_2<\nu_1\le1 \) (see e.g., \cite{MR1279172}). In particular, we can choose \(A=\Gamma\), where \(\Gamma\) is a Jordan closed curve.
\end{remark}

\subsection{Douglis algebras and hyperanalytic functions}

A Douglis algebra is an algebra of hypercomplex numbers generated by the elements \(i\) and \(e\); it is a generalization of the complex numbers in the complex plane \(\mathbb{C}\). For details about this topic we refer the reader to \cite{GilbertBuchanan1983, douglis1953, MR150320}. 
Let \(\mathbb{D}\) be the Douglis algebra; the multiplication in \(\mathbb{D}\) is governed by the rules:
\[
i^2 = -1,\quad ie = ei,\quad e^r = 0,\quad e^0 = 1,
\]
where \(r\) is a positive integer. Any arbitrary element \(a \in \mathbb{D}\) may be written as a hypercomplex number in the form
\[
a = \sum_{k=0}^{r-1} a_k e^k,
\]
where each \(a_k\) is a complex number; \(a_0\) is the complex part of \(a\), and \(A = \sum_{k=1}^{r-1} a_k e^k\) is its nilpotent part.

Let \(a\) be any hypercomplex number, then its conjugation \(\overline{a}\) is defined as \(\overline{a} = \sum_{k=0}^{r-1} \overline{a_k} e^k\). The algebraic norm in \(\mathbb{D}\) is defined by \(|a| := \sum_{k=0}^{r-1} |a_k|\).

If the complex part \(a_0\) of a hypercomplex number \(a\) is nonzero, then the multiplicative inverse \(a^{-1}\) of \(a\) is given by
\[
a^{-1} = a_0^{-1} \sum_{k=0}^{r-1} (-1)^k \left(\frac{A}{a_0}\right)^k.
\]
If \(a_0 = 0\) then \(a\) is called nilpotent and it does not have a multiplicative inverse.

Let \(f\) be a \(\mathbb{D}\)-valued function; then \(f\) may be written as \(f = \sum_{k=0}^{r-1} f_k e^k\), where \(f_k\) are complex-valued functions. Each time we assign a property such as continuity, differentiability, integrability, etc., to \(f\) it is meant that all components \(f_k\) share this property.

The Douglis operator \(\overline{\partial}_q\) is given by
\[
\overline{\partial}_q := \overline{\partial} + q(z)\partial, \quad z = x + iy,
\]
where \(q(z)\) is a known nilpotent hypercomplex function and
\[
\overline{\partial} := \frac{1}{2}(\partial_x + i\partial_y),\quad \partial := \frac{1}{2}(\partial_x - i\partial_y).
\]
\begin{definition}
Let \(\Omega\subset \mathbb{C}\) be a domain. A smooth hypercomplex function \(f\) defined in \(\Omega\) is said to be hyperanalytic in \(\Omega\) if \(\overline{\partial}_q f = 0\) in \(\Omega\). The class of hyperanalytic functions on \(\Omega\) is denoted by \(\mathcal{H}_q(\Omega)\).
\end{definition}

As an example of a hyperanalytic function we take the generating solution of the Douglis operator given by
\[
W(z) = z + \sum_{k=1}^{r-1} W_k(z) e^k,
\]
where its nilpotent part possesses bounded and continuous derivatives up to order two in \(\mathbb{C}\). The following properties concerning the generating solution will be used frequently:
\[
|W(z_1) - W(z_2)| \asymp |z_1 - z_2|,\qquad |W(z_1) - W(z_2)|^{-1} \leq c|z_1 - z_2|^{-1},\; z_1 \neq z_2,
\]
where \(c\) denotes a constant which may vary from one occurrence to the next; in general these constants only depend on \(q\).

Another important example of a hyperanalytic function is the so-called hypercomplex Cauchy kernel, i.e., the fundamental solution of the Douglis operator, given by
\[
e_z(\zeta) := \frac{1}{2\pi} \frac{\partial_\zeta W(\zeta)}{W(\zeta) - W(z)}, \quad \zeta \neq z.
\]
The nature of the singularity of \(e_z(\zeta)\) is the same as that which the complex Cauchy kernel \(\frac{1}{\zeta - z}\) has at \(\zeta = z\).

We will denote by \(C(\Gamma)\) the set of all continuous hypercomplex functions defined on \(\Gamma\). Moreover, for \(f \in C(\Gamma)\) we introduce the modulus of continuity for \(f\):
\[
\omega_f(\varepsilon) := \sup_{\substack{t_1,t_2 \in \Gamma, |t_1-t_2| \leq \varepsilon}} \max|f(t_1) - f(t_2)|, \quad \varepsilon > 0.
\]
Let us consider also the subclass \(H^\nu(\Gamma) \subset C(\Gamma)\) of all functions \(f\) satisfying a H\"{o}lder condition \(\omega_f(\varepsilon) \leq c \varepsilon^\nu\) with exponent \(\nu\), \(0 < \nu \leq 1\). Here and subsequently, \(\mathrm{diam}(E)\) denotes the diameter of \(E \subset \mathbb{C}\).

For a \(\mathbb{D}\)-valued function \(f = \sum_{k=0}^{r-1} f_k e^k \in H^\nu(\Gamma)\), we mean that each complex component \(f_k \in H^\nu(\Gamma)\).

It is worth pointing out that a very successful tool in the theory of Riemann boundary value problems both for analytic and hyperanalytic functions is the corresponding Cauchy type integral.

If \(\Gamma\) is a Jordan closed rectifiable curve, then for any \(f \in C(\Gamma)\), the customary hypercomplex Cauchy type integral
\begin{equation}\label{hypercomplex_Cauchy_rectifiable}
(\mathcal{C}_\Gamma f)(z) := \int_\Gamma e_z(\zeta) n_q(\zeta) f(\zeta) \, \mathrm{d}s, \quad z \notin \Gamma,
\end{equation}
where \(n_q(\zeta) := n(\zeta) + \overline{n(\zeta)}q(\zeta)\) with \(n(\zeta)\) being the exterior unit normal vector at the point \(\zeta\) on \(\Gamma\) in Federer's sense (see \cite{MR503901}), and \(\mathrm{d}s\) denotes the arclength differential, exists and represents a hyperanalytic function in \(\mathbb{C} \setminus \Gamma\).

At almost all (with respect to \(\mathrm{d}s\)) points \(t \in \Gamma\) this function has non-tangential boundary limit values from both sides, and these values almost everywhere satisfy the relation
\[
(\mathcal{C}_\Gamma f)^+(t) - (\mathcal{C}_\Gamma f)^-(t) = f(t), \quad t \in \Gamma.
\]

If \(f \in H^\nu(\Gamma)\) and \(\nu > 1/2\), then the function \(\mathcal{C}_\Gamma f\) has continuous boundary values on the whole \(\Gamma\) (see \cite{Zeng1990}). If the curve is Ahlfors-David regular, then these properties are valid for any \(\nu \in (0, 1]\) (see \cite{MR2005hyperRBVP}).

\subsection{Box dimension and $d$-summable sets in $\mathbb{C}$}

The standard approach considers a fractal to be a set with a non-integer Hausdorff dimension. However, frequently the box dimension is a more appropriate dimension than the Hausdorff dimension to measure the roughness of a bounded set.

\begin{definition}\label{box_dimension}
Let \( \Gamma \) be a Jordan closed curve, then
\[
\mathrm{Dmb}\,\Gamma = \limsup\limits_{\varepsilon\to 0}\frac{\log{N(\varepsilon,\Gamma)}}{-\log{\varepsilon}}
\]
is called the box dimension of \( \Gamma \), where \(N(\varepsilon,\Gamma)\) is the least number of disks with radius not exceeding \(\varepsilon\) needed to cover \(\Gamma\). Let \(\delta>0\), the family \(\{U_i\}_{i=0}^\infty\) is called a \(\delta\)-cover of \( \Gamma \), if \( \Gamma \subseteq \bigcup\limits_{i=0}^{\infty} U_i\) and \( \mathrm{diam}(U_i)=\sup\limits_{z_1,z_2\in U_i}|z_1 - z_2|\le \delta \) for each \( i \). Then
\[\mathrm{Dmh}\, \Gamma = \inf\{ s \geq 0 : \mathscr{H}^s(\Gamma) = 0 \} = \sup\{ s : \mathscr{H}^s(\Gamma) = \infty \}\]
is called the Hausdorff dimension of \( \Gamma \), where
\[
\mathscr{H}^s(\Gamma) = \lim\limits_{\delta \to 0} \mathscr{H}_\delta^s(\Gamma)=\lim\limits_{\delta \to 0}\inf\left\{ \sum_{i=0}^\infty [\mathrm{diam}(U_i)]^s : \{U_i\}_{i=0}^\infty \text{ is a } \delta\text{-cover of } \Gamma \right\}
\]
is the Hausdorff measure of \( \Gamma \) (see e.g., \cite{MR1102677}).
\end{definition}
\begin{remark}\label{properties_of_dimension}
For the Jordan closed curve \( \Gamma \), we have that \( 1\le\mathrm{Dmh}\, \Gamma\le\mathrm{Dmb}\, \Gamma\le2 \). And if \( \Gamma \) is a rectifiable Jordan closed curve, then \( \mathrm{Dmb} \, \Gamma =1 \) (see e.g., \cite{Kat83}).
\end{remark}

\begin{definition}
The closed Jordan curve \( \Gamma \) is said to be \( d \)-summable if the improper integral
\[
\int_0^1 N(x,\Gamma) x^{d-1} \mathrm{d}x
\]
converges, where \(N(x,\Gamma)\) is defined as in Definition \ref{box_dimension}.
\end{definition}
\begin{remark}

This notion was introduced by Jenny Harrison and Alec Norton in \cite{greenthm, Geometintegrat}. It is a geometric property describing a set (such as a fractal boundary) for which the sum of the diameters of its covering sets scales in a specific way relative to a dimension $d$. Moreover, if \(\Gamma\) is \(d\)-summable, then \(\mathrm{Dmb}\, \Gamma \le d\), and if \(\mathrm{Dmb}\, \Gamma < d\), then \(\Gamma\) is \(d\)-summable. Thus, every rectifiable Jordan closed curve is \(d\)-summable for every \(d>1\), and a fractal curve is a \(d\)-summable curve for some \(d\). This concept can be generalized to \(h\)-summability. For thorough treatment we refer the reader to \cite{MR3146343}.
\end{remark}

\subsection{Whitney extension}

One of the techniques for solving Riemann boundary value problems on non-rectifiable curves is to extend the domain of definition of \( f \) to the entire complex plane. This technique is called Whitney extension (see e.g., \cite{MR290095}). Let \( \Gamma \) be a Jordan closed curve. Then there exists a collection of closed cubes \( \mathcal{W} \), \( \mathcal{W} = \{ Q_1, Q_2, \ldots, Q_k, \cdots \} \) such that
\begin{itemize}
\item \( \displaystyle \bigcup_{k=1}^{\infty} Q_k = \Omega = \mathbb{C}\setminus\Gamma \),
\item The interiors of \( Q_k \) are mutually disjoint,
\item \( c_1 \mathrm{diam}(Q_k) \leq \mathrm{dist}(Q_k, \Gamma) \leq c_2 \mathrm{diam}(Q_k) \),
\end{itemize}
where \( \mathrm{dist}(Q_k, \Gamma) = \inf\limits_{z\in Q_k,t\in \Gamma}|z-t| \).
The constants \( c_1 \) and \( c_2 \) are independent of \( \Gamma \). In fact we may take \( c_1 = 1 \) and \( c_2 = 4 \). Now let \( Q_0 \) denote the cube of unit length centered at the origin. Fix \(\varepsilon>0\) and a \( C^\infty \) function \( \varphi \) with the properties: \( 0 \leq \varphi \leq 1 \); \( \varphi(z) = 1 \) for \( z \in Q_0 \); and \( \varphi(z) = 0 \) for \( z \notin (1 + \varepsilon)Q_0 \). Let
\[
\varphi_k(z) = \varphi\left( \frac{z - z^k}{l_k} \right),
\]
where \( z^k \) is the center of \( Q_k \) and \( l_k \) is the common length of its sides. We now define \( \varphi_k^*(z) \) for \( z \in \mathbb{C}\setminus\Gamma \) by
\[
\varphi_k^*(z) = \frac{\varphi_k(z)}{\Phi(z)},
\]
where \( \Phi(z) = \sum\limits_{k=1}^\infty \varphi_k(z)\).
\begin{remark}
In the expression \( \Phi(z) = \sum\limits_{k=1}^\infty \varphi_k(z) \), the summation is always a finite sum, i.e., only finitely many terms are non-zero in the sum (see e.g., \cite{MR290095}). The obvious identity
\[
\sum\limits_{k=1}^\infty \varphi_k^*(z) \equiv 1, \quad z \in \mathbb{C}\setminus\Gamma,
\]
which gives a partition of unity.
\end{remark}

For the given function \( f \in H^\nu(\Gamma) \) (in the hypercomplex sense, i.e., each component is in \(H^\nu(\Gamma)\)), consider the function \( \mathcal{E}_0(f) \) defined by
\begin{equation}\nonumber
\mathcal{E}_0(f)(z) =
\begin{cases}
f(z), &z \in \Gamma, \\
\sum\limits_{k=1}^\infty f(p_k) \varphi_k^*(z), &z \in \mathbb{C}\setminus\Gamma,
\end{cases}
\end{equation}
where \( p_k \) satisfies the property that \( \mathrm{dist}(Q_k, \Gamma) = \mathrm{dist}(Q_k, p_k) \). \( \mathcal{E}_0(f) \) is called the Whitney extension of \( f \). Now let
\[
  f^\omega(z)=\chi(z)\mathcal{E}_0(f)(z),
\]
where \( \chi(z) \) is the characteristic function of \( \overline{D^+} \), i.e.,
\[
\chi(z)=
\begin{cases}
1, &z \in \overline{D^+},\\
0, &z \in D^-.
\end{cases}
\]
Then we have that \( f^\omega \in H^\nu(\overline{D^+}) \), \( f^\omega \) has partial derivatives with respect to \( x \) and \( y \) on \( D^+ \) and
\begin{equation}\label{up_bounded_derivative}
 \left| \overline{\partial}_q f^\omega(z) \right| \leq C [\mathrm{dist}(z,\Gamma)]^{\nu - 1}, \quad \forall z \in D^+,
\end{equation}
where \( C \) is a constant (see e.g., \cite{MR290095, Aa32}).
\begin{remark}
We utilize the Whitney decomposition, enabling functions defined on the curve to be extended to the entire complex plane. Consequently, the originally undefined curvilinear integrals can be converted into well-defined domain integrals. This also provides a method for handling fractal problems.
\end{remark}

\subsection{The hypercomplex Cauchy type integral on $d$-summable curves}

When assuming a much more pathological situation, e.g., \(\Gamma\) is assumed to be a fractal, then the definition (\ref{hypercomplex_Cauchy_rectifiable}) of the Cauchy type integral fails. We now present an alternative definition of the hypercomplex Cauchy type integral when the contour is allowed to be a fractal, following \cite{Aa32}.

\begin{definition}\label{hypercomplex_Cauchy_dsummable}
Let \( d \in (1,2) \) and let \(\Gamma\) be a \( d \)-summable Jordan closed curve, and suppose \( \nu > d - 1 \). The hypercomplex Cauchy type integral of \( f \in H^\nu(\Gamma) \) is defined by the formula
\begin{equation}\label{Cauchy_hypercomplex_integral}
\mathcal{C}_\Gamma (f)(z) = f^\omega(z) - \iint_{D^+} e_z(\zeta) \overline{\partial}_q \mathcal{E}_0(f)(\zeta) \, \mathrm{d}\xi \mathrm{d}\eta, \quad z \in \mathbb{C} \setminus \Gamma,
\end{equation}
with \( \zeta = \xi + i\eta \).
\end{definition}

The following proposition makes this definition legitimate (see \cite{Aa32}).

\begin{proposition}\label{prop_Cauchy_well_defined}
The hypercomplex function defined by (\ref{Cauchy_hypercomplex_integral}) is correctly defined for any \( z \in \mathbb{C} \setminus \Gamma \) and its value does not depend on the particular choice of \( \mathcal{E}_0(f) \).
\end{proposition}

We now state two key theorems from \cite{Aa32} concerning the jump problem.

\begin{theorem}\label{thm_jump_continuous}
Let \( \Gamma \) be a \( d \)-summable curve and \( f \in H^\nu(\Gamma) \). If \( \nu > \frac{d}{2} \), then the Cauchy type integral (\ref{Cauchy_hypercomplex_integral}) has continuous limit values on \( \Gamma \) from both domains \( D^\pm \). Moreover,
\[
(\mathcal{C}_\Gamma f)^+(t) - (\mathcal{C}_\Gamma f)^-(t) = f(t), \quad t \in \Gamma.
\]
\end{theorem}

To ensure uniqueness of the solution of the jump problem we need to introduce some additional requirements. A function \(\Phi\), hyperanalytic in \(\mathbb{C}\setminus\Gamma\), must satisfy a H\"{o}lder condition with exponent \(\mu\), \(0 < \mu < 1\) on each of the sets \(\overline{D^\pm}\), i.e., the restrictions \(\Phi|_{D^+}\) and \(\Phi|_{D^-}\) must be \(\mu\)-H\"{o}lder continuous in the closed domains respectively, and the boundary values of these restrictions \(\Phi^\pm\) are the usual continuous limit values. By a Dolzhenko-type theorem for hyperanalytic functions (see \cite{Aa32}), the desired uniqueness follows from the removability of the curve \(\Gamma\) under the condition
\begin{equation}\label{dolzhenko_condition}
\mu > \mathrm{Dmh}\,\Gamma - 1.
\end{equation}

It is essential to point out that the integral term in (\ref{Cauchy_hypercomplex_integral}) belongs to \(H^\mu(\mathbb{C})\) with
\begin{equation}\label{holder_exponent_bound}
\mu < \frac{2\nu - d}{2 - d}.
\end{equation}

A function \(\Phi\), being \(\mu\)-H\"{o}lder continuous on \(\overline{D^\pm}\) whenever \(\mu\) satisfies (\ref{dolzhenko_condition}) and (\ref{holder_exponent_bound}), is said to be of class \(\mathcal{H}_\mu\).

\begin{theorem}\label{thm_jump_unique}
Under the hypotheses of Theorem \ref{thm_jump_continuous}, if moreover
\[
\mathrm{Dmh}\,\Gamma - 1 < \mu < \frac{2\nu - d}{2 - d},
\]
then there exists a unique solution of the jump problem
\[
\Phi^+(t) - \Phi^-(t) = f(t), \quad t \in \Gamma,
\]
in the class \(\mathcal{H}_\mu\).
\end{theorem}

\subsection{The hyperanalytic Riemann boundary value problem}

Let \(\Gamma\) be a Jordan closed curve. \(F \in \mathcal{H}_q(\mathbb{C}\setminus\Gamma)\) is called sectionally hyperanalytic function with the jump \(\Gamma\), if the boundary values \(F^\pm(t)\) exist for each \(t \in \Gamma\).

\begin{problem}[Riemann boundary value problem for hyperanalytic functions on a $d$-summable Jordan closed curve]\label{hyper_RBVP}
Let \( \Gamma \) be a $d$-summable Jordan closed curve. Find a sectionally hyperanalytic function \( F \) with the jump \( \Gamma \) such that
\begin{equation}\label{RBV_equation}
F^+(t)=G(t)F^-(t)+f(t), \quad \forall t\in\Gamma,
\end{equation}
where \( G,f \in H^\nu(\Gamma) \) are given for \( \nu\in (0,1] \), \(G\) being a hypercomplex function whose complex part \(G_0\) never vanishes on \(\Gamma\), and the operations in (\ref{RBV_equation}) are in the sense of the Douglis algebra.
\end{problem}

Let \( \kappa = \frac{1}{2\pi}\Delta_\Gamma \arg(G_0(t)) \), where \(G_0\) is the complex part of \(G\). For \( \mu\in(0,1] \), define the functions class
\begin{equation}\nonumber
\begin{split}
\mathcal{R}_m^\mu(\overline{\partial}_q) =\{F \in H^\mu(\overline{D^+})\cap H^\mu(\overline{D^-}): &F \text{ is a sectionally hyperanalytic function with jump }\Gamma, \\
&\text{and }\mathrm{Ord}(F,\infty)\le m\}.
\end{split}
\end{equation}
Let \( \mathscr{P}_\kappa(W) \) denote the set of all hypercomplex polynomials in \(W(z)\) of degree not exceeding \( \kappa \), i.e., polynomials of the form \(\sum_{j=0}^{\kappa} a_j W^j(z)\) with \(a_j \in \mathbb{D}\).

The following theorem from \cite{Aa32} solves Problem \ref{hyper_RBVP}.

\begin{theorem}[Abreu Blaya, Bory Reyes, and Vilaire]\label{hyper_theorem}
Suppose that \( d\in(1,2) \). If \( \nu > \dfrac{d}{2} \) and \( \mathrm{Dmh}\, \Gamma - 1 < \mu < \dfrac{2\nu - d}{2 - d} \) or \( \nu = 1 \) and \( \mathrm{Dmh}\, \Gamma - 1 < \mu < 1 \), then, when \( \kappa \ge 0 \), the general solution of Problem \textup{\ref{hyper_RBVP}} in the class \( \mathcal{R}_0^\mu(\overline{\partial}_q) \) is of the form \( F(z) = \Phi(z) + X(z)P_\kappa(W(z)) \), where \(P_\kappa\in\mathscr{P}_\kappa(W)\). When \( \kappa = -1 \), \( \Phi(z) \) is the unique solution of Problem \textup{\ref{hyper_RBVP}} in the class \( \mathcal{R}_0^\mu(\overline{\partial}_q) \). Finally, when \( \kappa < -1 \), the necessary and sufficient condition for Problem \textup{\ref{hyper_RBVP}} to be solvable in the class \(\mathcal{R}_0^\mu(\overline{\partial}_q) \) is
\begin{equation}\label{hyper_solvable_condition}
    \iint_{D^+} \frac{\overline{\partial}_q f^\omega(\zeta)}{X^+(\zeta)} W^{m-1}(\zeta) \, \mathrm{d}\xi \mathrm{d}\eta = 0, \quad m = 1, 2, \cdots, -\kappa - 1,
\end{equation}
when condition \((\ref{hyper_solvable_condition})\) is satisfied, the function \( \Phi(z) \) is also the unique solution of Problem \textup{\ref{hyper_RBVP}} in the class \( \mathcal{R}_0^\mu(\overline{\partial}_q) \). The explicit expression of \( \Phi(z) \) is
\begin{equation}\label{hyper_Phi}
    \Phi(z) = f^\omega(z) + X(z) \iint_{D^+} e_z(\zeta) \left( \frac{\overline{\partial}_q f^\omega}{X} \right)(\zeta) \, \mathrm{d}\xi \mathrm{d}\eta,
\end{equation}
where
\begin{equation}\label{hyper_X}
X(z) =
\begin{cases}
\exp\{ \mathcal{C}_\Gamma (\ln[W^{-\kappa} G])(z) \}, & z\in D^+, \\
W(z)^{-\kappa} \exp\{ \mathcal{C}_\Gamma (\ln[W^{-\kappa} G])(z) \}, & z\in D^-,
\end{cases}
\end{equation}
and the hypercomplex exponential and logarithmic functions are understood in the sense of the Douglis algebra (see \cite{GilbertBuchanan1983}). Here \(\ln[W^{-\kappa} G]\) denotes a continuous branch of the hypercomplex logarithm such that its complex part is single-valued.
\end{theorem}

\section{Riemann boundary value problem for poly-hyperanalytic functions}

In this section, we mainly solve the Riemann boundary value problem for poly-hyperanalytic functions. First, we give the definition of poly-hyperanalytic functions.

\begin{definition}
Let \( \Omega \) be an open set in \( \mathbb{C} \), \( F: \Omega \to \mathbb{D} \) is called poly-hyperanalytic function of order \( n \), if \( F \in C^n(\Omega) \) and
\[ \overline{\partial}_q^{\,n}F(z) = 0 , \quad \forall z \in \Omega. \]
The class of poly-hyperanalytic functions of order \( n \) on \( \Omega \) is denoted by \( \mathcal{H}_q^n(\Omega) \).
Particularly, \( \mathcal{H}_q^1(\Omega)=\mathcal{H}_q(\Omega) \) is the set of hyperanalytic functions on \( \Omega \).
\end{definition}

Let \( \Gamma \) be a Jordan closed curve. \( F \) is called a sectionally poly-hyperanalytic function of order \( n \) with the jump \( \Gamma \), if \( F \in \mathcal{H}_q^n(\mathbb{C} \setminus \Gamma) \) and \( (\overline{\partial}_q^{\,k} F)^\pm(t) \) exist for all \( t\in\Gamma \) and \( k=0,1, \cdots,n-1 \).

It is easy to see that
\[
\overline{\partial}_q\left[ F(z)G(z) \right] = \overline{\partial}\left[F(z)G(z)\right]+q(z)\partial\left[F(z)G(z)\right]=G(z)\overline{\partial}_q F(z)+F(z)\overline{\partial}_q G(z), \quad z\in \mathbb{C}\setminus\Gamma,
\]
for \(F,G\in C^n(\mathbb{C}\setminus\Gamma)\) by the Leibniz rule for \(\overline{\partial},\partial\). 
Moreover,
\[ \overline{\partial}_q^{\,n}\left[ F(z)G(z) \right] = F(z)\overline{\partial}_q^{\,n}G(z), \quad z\in \mathbb{C}\setminus\Gamma \]
for \( F \in \mathcal{H}_q(\mathbb{C}\setminus\Gamma) \) and \( G\in C^n(\mathbb{C}\setminus\Gamma) \).

Now observe that for the conjugate variable \( \overline{z} = x - iy \), we have \( \partial(\overline{z}^k) = 0 \) for any integer \(k\). Consequently,
\[
\overline{\partial}_q (\overline{z}^k) = \overline{\partial}(\overline{z}^k) + q(z)\partial(\overline{z}^k) = k\overline{z}^{k-1}.
\]
Then one can see that \( \sum\limits_{k = 0}^{n - 1}\overline{z}^k F_k(z)\) \(\in\) \(\mathcal{H}_q^n(\mathbb{C}\setminus\Gamma) \), where \( F_k \in \mathcal{H}_q(\mathbb{C}\setminus\Gamma) \) for \( k = 0, 1, \cdots, n - 1 \). Actually, we have the following decomposition result, which is analogous to the decomposition for poly-analytic and poly-\(\beta\)-analytic functions.

\begin{theorem}\label{poly_hyper_decomposition}
  \( \mathcal{H}_q^n(\mathbb{C}\setminus\Gamma) = \mathcal{H}_q(\mathbb{C}\setminus\Gamma) \oplus \overline{z}\mathcal{H}_q(\mathbb{C}\setminus\Gamma) \oplus \cdots \oplus \overline{z}^{n - 1}\mathcal{H}_q(\mathbb{C}\setminus\Gamma) \).
\end{theorem}
\begin{Proof}
We have shown that \( \sum\limits_{k = 0}^{n - 1}\overline{z}^k F_k(z) \in \mathcal{H}_q^n(\mathbb{C}\setminus\Gamma) \) if \( F_k \in \mathcal{H}_q(\mathbb{C}\setminus\Gamma) \) for \( k = 0, 1, \cdots, n - 1 \). Suppose that for each \( F\in \mathcal{H}_q^{n-1}(\mathbb{C}\setminus\Gamma) \), there exists \( F_k \in \mathcal{H}_q(\mathbb{C}\setminus\Gamma) \) for \( k = 0, 1, \cdots, n - 2 \) such that \( F(z) = \sum\limits_{k = 0}^{n - 2}\overline{z}^k F_k(z) \). Let \(F\in\mathcal{H}_q^{n}(\mathbb{C}\setminus\Gamma)\), \( E(z) = \frac{1}{(n - 1)!}\overline{\partial}_q^{\,n-1} F(z) \), and \( G(z) = F(z) - \overline{z}^{n - 1}E(z) \), then we have \( \overline{\partial}_q^{\,n-1}G(z) = 0 \), which means that
\[ G(z) = \sum\limits_{k = 0}^{n - 2}\overline{z}^k G_k(z) \text{ and }  F(z) = \sum\limits_{k = 0}^{n - 2}\overline{z}^k G_k(z) + \overline{z}^{n - 1}E(z),\]
where \( G_k, E \in \mathcal{H}_q(\mathbb{C}\setminus\Gamma) \) and \( k = 0, 1, \cdots, n-2 \). By induction, we show that \[ \mathcal{H}_q^n(\mathbb{C}\setminus\Gamma) = \mathcal{H}_q(\mathbb{C}\setminus\Gamma) + \overline{z}\mathcal{H}_q(\mathbb{C}\setminus\Gamma) + \cdots + \overline{z}^{n - 1}\mathcal{H}_q(\mathbb{C}\setminus\Gamma). \]
And it is worth pointing out that \(\mathcal{H}_q(\mathbb{C}\setminus\Gamma), \overline{z}\mathcal{H}_q(\mathbb{C}\setminus\Gamma),...,\overline{z}^{n-1}\mathcal{H}_q(\mathbb{C}\setminus\Gamma)\) are linear spaces over \(\mathbb{D}\). Suppose that \( F(z) = \sum\limits_{k = 0}^{n - 1}\overline{z}^k F_k(z) = 0 \), \( F_k \in \mathcal{H}_q(\mathbb{C}\setminus\Gamma) \), \( k = 0, 1, \cdots, n - 1 \), then \( F_{n - 1}(z) = 0 \) by \( \overline{\partial}_q^{\,n-1}F(z) = 0 \) and \( F_{n - 2}(z) = 0 \) by \( \overline{\partial}_q^{\,n-2}F(z) = 0 \) and \( \cdots \cdots \). Hence for each \( k \), \( F_k(z) = 0 \), which means that \[ \mathcal{H}_q^n(\mathbb{C}\setminus\Gamma) = \mathcal{H}_q(\mathbb{C}\setminus\Gamma) \oplus \overline{z}\mathcal{H}_q(\mathbb{C}\setminus\Gamma) \oplus \cdots \oplus \overline{z}^{n - 1}\mathcal{H}_q(\mathbb{C}\setminus\Gamma). \]
\end{Proof}

The specific statement of the Riemann boundary value problem for poly-hyperanalytic functions is as follows.

\begin{problem}[Riemann boundary value problem for poly-hyperanalytic functions on a Jordan closed curve]\label{Jordan_curve_poly_hyper_RBVP}
Let \( \Gamma \) be a Jordan closed curve. Find a sectionally poly-hyperanalytic function \( F \) with the jump \( \Gamma \) such that
\begin{equation}\label{poly-hyperanalytic_RQ}
\begin{cases}
\begin{aligned}
F^+(t) &= G(t)F^-(t) + f_0(t), \\
(\overline{\partial}_q F)^+(t) &= G(t)(\overline{\partial}_q F)^-(t) + f_1(t), \\
\vdots & \\
(\overline{\partial}_q^{\,n-1} F)^+(t) &= G(t)(\overline{\partial}_q^{\,n-1} F)^-(t) + f_{n-1}(t),
\end{aligned}
\ t \in \Gamma,
\end{cases}
\end{equation}
where \( G,f_k \in H^\nu(\Gamma) \) are given for \( \nu\in (0,1] \), \( k=0,1,\cdots,n-1 \), \(G\) is a hypercomplex function whose complex part \(G_0\) never vanishes on \(\Gamma\).
\end{problem}

Next, to address this problem, we introduce a class of functions. For \( \mu\in(0,1] \), define the functions class
\begin{equation}\nonumber
\begin{split}
\mathcal{R}_m^\mu(\overline{\partial}_q^{\,n}) =\{ F : &  F(z)=\sum\limits_{k=0}^{n-1}\overline{z}^kF_k(z) \text{ is a sectionally poly-hyperanalytic function with jump }\Gamma, \\
&\mathrm{Ord}(F_k,\infty)\le m \text{ and } F_k \in H^\mu(\overline{D^+})\cap H^\mu(\overline{D^-}) \text{ for } k=0,1,\cdots,n-1\}.
\end{split}
\end{equation}

\begin{lemma}\label{lemma_of_poly_hyper_class_equivalent}
The necessary and sufficient condition for \( F(z)=\sum\limits_{k=0}^{n-1}\overline{z}^kF_k(z)\in \mathcal{R}_0^\mu(\overline{\partial}_q^{\,n}) \) and \( F \) satisfying the system of equations \((\ref{poly-hyperanalytic_RQ})\) is that for each \( k=0,1,\cdots,n-1 \), \( F_k \in \mathcal{R}_0^\mu(\overline{\partial}_q) \) and \( F_k^+(t)=G(t)F_k^-(t)+\sum\limits_{j=k}^{n-1}\dfrac{(-1)^{k+j}}{(j - k)!k!}\overline{t}^{j - k}f_j(t) \).
\end{lemma}
\begin{Proof}
Using Theorem \ref{poly_hyper_decomposition} we have that \( F(z)=\sum\limits_{k=0}^{n-1}\overline{z}^kF_k(z) \) where \( F_k \in \mathcal{H}_q(\mathbb{C}\setminus\Gamma) \) for \( k=0,1,\cdots,n-1 \), and \( \overline{\partial}_q^{\,k} F(z) = \sum\limits_{j=k}^{n-1} \frac{j!}{(j-k)!}\overline{z}^{j-k}F_j(z) \) for \( k=0,1,\cdots,n-1 \). Consequently, we obtain that \( F_k^\pm(t) \) exist for all \( t\in\Gamma \) and \(k=0,1,\cdots,n-1\). Conversely, for each \( k=0,1,\cdots,n-1 \), if \( F_k(z) \) is a sectionally hyperanalytic function with jump \( \Gamma \), then \( F(z)=\sum\limits_{k=0}^{n-1}\overline{z}^kF_k(z) \in \mathcal{H}_q^n(\mathbb{C}\setminus\Gamma) \) and \( (\overline{\partial}_q^{\,k} F)^\pm(t) \) exist for \( t\in\Gamma \), i.e., \(F\) is a sectionally poly-hyperanalytic function of order \(n\) with jump \( \Gamma \). Thus, by the definition, it is easy to see that \( F(z)=\sum\limits_{k=0}^{n-1}\overline{z}^kF_k(z)\in  \mathcal{R}_0^\mu(\overline{\partial}_q^{\,n}) \Longleftrightarrow F_k \in \mathcal{R}_0^\mu(\overline{\partial}_q) \) for any \( k=0,1,\cdots,n-1 \).

Now, let us consider (\ref{poly-hyperanalytic_RQ}). Since \( F(z)=\sum\limits_{k=0}^{n-1} \overline{z}^kF_k(z) \) where \( F_k\in \mathcal{H}_q(\mathbb{C}\setminus\Gamma) \) for \( k=0,1,\cdots,n-1 \), (\ref{poly-hyperanalytic_RQ}) is equivalent to
\begin{equation}\nonumber
\begin{cases}
\begin{aligned}
F_0^+(t) + \overline{t}F_1^+(t) + \cdots + \overline{t}^{n-1}F_{n-1}^+(t) &= G(t)(F_0^-(t) + \overline{t}F_1^-(t) + \cdots + \overline{t}^{n-1}F_{n-1}^-(t)) + f_0(t), \\
F_1^+(t) + \cdots + (n-1)\overline{t}^{n-2}F_{n-1}^+(t) &= G(t)( F_1^-(t) + \cdots + (n-1)\overline{t}^{n-2}F_{n-1}^-(t)) + f_1(t), \\
\vdots &\\
(n-1)!F_{n-1}^+(t) &= G(t)(n-1)!F_{n-1}^-(t) + f_{n-1}(t).
\end{aligned}
\end{cases}
\end{equation}
Let
\[
\mathbf{F}(z)=\begin{pmatrix}F_0(z)\\ F_1(z)\\ \vdots \\ F_{n-1}(z)\end{pmatrix}, \quad z\in\mathbb{C}\setminus\Gamma,
\quad
\mathbf{G}(t)=\begin{pmatrix}
G(t) & \cdots & 0 \\
\vdots & \ddots & \vdots \\
0 & \cdots & G(t)
\end{pmatrix},
\quad
\mathbf{f}(t)=\begin{pmatrix}f_0(t)\\ f_1(t)\\\vdots\\ f_{n-1}(t)\end{pmatrix}, \quad t\in \Gamma.
\]
Then the equations above are equivalent to the following matrix equation
\begin{equation}\label{poly-hyperanalytic_matrixequation_A}
  \mathbf{A}(t)\mathbf{F}^+(t) = \mathbf{A}(t)\mathbf{G}(t)\mathbf{F}^-(t) + \mathbf{f}(t), \quad t\in \Gamma,
\end{equation}
where
\begin{equation}\label{defination_of_A_hyper}
\mathbf{A}(z) = \left(a_{ij}(z)\right)_{0\le i,j\le n-1} \text{ with } a_{ij}(z) =
\begin{cases}
\dfrac{j!}{(j-i)!} \overline{z}^{j - i}, & i \le j, \\
0, & i > j,
\end{cases}
 \quad z\in\mathbb{C}.
\end{equation}
It is easy to see that \(\mathbf{A}^{-1}(z) = \mathbf{B}(z)\) where
\[
\mathbf{B}(z) = \left(b_{ij}(z)\right)_{0\le i,j\le n-1} \text{ with } b_{ij}(z) =
\begin{cases}
\dfrac{(-1)^{i+j}}{(j - i)!i!}\overline{z}^{j - i}, & i \le j, \\
0, & i > j,
\end{cases} \quad z\in\mathbb{C}.
\]
Therefore, (\ref{poly-hyperanalytic_matrixequation_A}) is equivalent to the following matrix equation
\begin{equation}\nonumber
  \mathbf{F}^+(t)=\mathbf{G}(t)\mathbf{F}^-(t)+\mathbf{B}(t)\mathbf{f}(t), \quad t\in\Gamma.
\end{equation}
Let
\begin{equation}\label{checkfk_hyper}
\check{f}_k(t)=\sum_{j=k}^{n-1} b_{kj}(t)f_j(t)=\sum_{j=k}^{n-1}\dfrac{(-1)^{k+j}}{(j - k)!k!}\overline{t}^{j - k}f_j(t), \quad k=0,1,\cdots,n-1.
\end{equation}
Thus, we obtain that
\begin{equation}\label{poly-hyperanalytic_componet_equation}
\begin{cases}
\begin{aligned}
F_0^+(t) &= G(t)F_0^-(t) + \check{f}_0(t), \\
F_1^+(t) &= G(t)F_1^-(t) + \check{f}_1(t), \\
&\vdots \\
F_{n-1}^+(t) &= G(t)F_{n-1}^-(t) + \check{f}_{n-1}(t),
\end{aligned}
\end{cases} \ t \in \Gamma,
\end{equation}
which are \(n\) independent Riemann boundary value problems for hyperanalytic functions. In the above discussion, we have shown that \(F\) satisfies (\ref{poly-hyperanalytic_RQ}) if and only if for \( k=0,1,\cdots,n-1 \), \( F_k \) satisfies \( F_k^+(t)=G(t)F_k^-(t)+\check{f}_k(t) \). Thus, the lemma holds.
\end{Proof}

Now, we present the solution to Problem \ref{Jordan_curve_poly_hyper_RBVP}. Let \( \kappa = \frac{1}{2\pi}\Delta_\Gamma \arg(G_0(t)) \), where \(G_0\) is the complex part of \(G\).

\begin{theorem}\label{solution_of_poly_hyper_RBVP}
Suppose that \(\Gamma\) is \(d\)-summable and \(d\in(1,2)\). If \( \nu > \dfrac{d}{2} \) and \( \mathrm{Dmh}\, \Gamma - 1 < \mu < \dfrac{2\nu - d}{2 - d} \) or \( \nu = 1 \) and \( \mathrm{Dmh}\, \Gamma - 1 < \mu < 1 \), then, when \( \kappa \ge 0 \), the general solution of Problem \textup{\ref{Jordan_curve_poly_hyper_RBVP}} in the class \( \mathcal{R}_0^\mu(\overline{\partial}_q^{\,n}) \) is of the form \( F(z) = \sum\limits_{j=0}^{n-1}\overline{z}^j\Phi_j(z) + X(z)Q_{n-1,\kappa}(z,\overline{z}) \), where
\begin{equation}\label{Phi_j_hyper}
    \Phi_j(z) = \check{f}^\omega_j(z) + X(z) \iint_{D^+} e_z(\zeta) \left( \frac{\overline{\partial}_q \check{f}^\omega_j}{X} \right)(\zeta) \, \mathrm{d}\xi \mathrm{d}\eta, \quad j=0,1,\cdots,n-1,
\end{equation}
\(\check{f}_j\) is defined by \((\ref{checkfk_hyper})\), \( X \) is defined as in \((\ref{hyper_X})\), and \( Q_{n-1,\kappa}(z,\overline{z}) \in \mathscr{P}_\kappa(W)+\overline{z}\mathscr{P}_\kappa(W)+\cdots+\overline{z}^{n-1}\mathscr{P}_\kappa(W)\). When \( \kappa = -1 \), \( F(z) = \sum\limits_{j=0}^{n-1}\overline{z}^j\Phi_j(z) \) is the unique solution of Problem \textup{\ref{Jordan_curve_poly_hyper_RBVP}} in the class \( \mathcal{R}_0^\mu(\overline{\partial}_q^{\,n}) \). Finally, when \( \kappa < -1 \), the necessary and sufficient condition for Problem \textup{\ref{Jordan_curve_poly_hyper_RBVP}} to be solvable in the class \( \mathcal{R}_0^\mu(\overline{\partial}_q^{\,n}) \) is
\begin{equation}\label{poly_hyper_solvable_condition}
\iint_{D^+} \frac{\overline{\partial}_q \check{f}^\omega_j(\zeta)}{X^+(\zeta)} W^{m-1}(\zeta) \, \mathrm{d}\xi \mathrm{d}\eta = 0, \quad m = 1, 2, \cdots, -\kappa - 1,\ j=0,1,\cdots,n-1,
\end{equation}
when condition \((\ref{poly_hyper_solvable_condition})\) is satisfied, the function \( F(z) = \sum\limits_{j=0}^{n-1}\overline{z}^j\Phi_j(z) \) is also the unique solution of Problem \textup{\ref{Jordan_curve_poly_hyper_RBVP}} in the class \( \mathcal{R}_0^\mu(\overline{\partial}_q^{\,n}) \).
\end{theorem}

\begin{Proof}
For \( j=0,1,\cdots,n-1 \), consider the Riemann boundary value problem \( F_j^+(t)=G(t)F_j^-(t)+\check{f}_j(t) \) and find solution \( F_j \) in \( \mathcal{R}_0^\mu(\overline{\partial}_q) \), where \( \check{f}_j \) is defined by (\ref{checkfk_hyper}). Moreover, by Remark \ref{properties_of_Holder_class}, \( \check{f}_j\in H^\nu(\Gamma) \) for \( j=0,1,\cdots,n-1 \). Note that since all operations are in the Douglis algebra and \(G(t)\) commutes with the scalar complex-valued coefficients in \(\mathbf{B}(t)\), the transformed boundary data \(\check{f}_j\) are well-defined hypercomplex functions.

Case 1, \( \kappa \ge 0 \). Using Theorem \ref{hyper_theorem}, we have that the general solution of \( F_j^+(t)=G(t)F_j^-(t)+\check{f}_j(t) \) in class \( \mathcal{R}_0^\mu(\overline{\partial}_q) \) is \[ F_j(z)=\Phi_j(z)+X(z)P_{j,\kappa}(W(z)),\quad j=0,1,\cdots,n-1, \]
where \(P_{j,\kappa}\in\mathscr{P}_\kappa(W)\), \( \Phi_j \) is defined by (\ref{Phi_j_hyper}) for \(j=0,1,\cdots,n-1\), and
\begin{equation}\nonumber
    X(z) =
    \begin{cases}
    \exp\{ \mathcal{C}_\Gamma (\ln[W^{-\kappa} G])(z) \}, & \,z\in D^+, \\
    W(z)^{-\kappa} \exp\{ \mathcal{C}_\Gamma (\ln[W^{-\kappa} G])(z) \}, & \,z\in D^-.
    \end{cases}
\end{equation}
Hence, by Lemma \ref{lemma_of_poly_hyper_class_equivalent}, we get the solution \( F \) of (\ref{poly-hyperanalytic_RQ}) in class \(\mathcal{R}_0^\mu(\overline{\partial}_q^{\,n})\) and
\begin{equation}\nonumber
F(z)=\sum_{j = 0}^{n - 1}\overline{z}^jF_j(z)=\sum\limits_{j = 0}^{n - 1}\overline{z}^j\Phi_j(z)+X(z)\sum\limits_{j = 0}^{n - 1}\overline{z}^jP_{j,\kappa}(W(z))=\sum\limits_{j = 0}^{n - 1}\overline{z}^j\Phi_j(z)+X(z)Q_{n - 1,\kappa}(z,\overline{z}),
\end{equation}
where \(Q_{n-1,\kappa}(z,\overline{z}) \in \mathscr{P}_\kappa(W)+\overline{z}\mathscr{P}_\kappa(W)+\cdots+\overline{z}^{n-1}\mathscr{P}_\kappa(W)\).

Case 2, \( \kappa=-1 \). Using Theorem \ref{hyper_theorem}, the equation \( F_j^+(t)=G(t)F_j^-(t)+\check{f}_j(t) \) has the unique solution \( F_j(z)=\Phi_j(z) \) in class \(\mathcal{R}_0^\mu(\overline{\partial}_q)\), \(j=0,1,\cdots,n-1\), where \( \Phi_j \) is defined by (\ref{Phi_j_hyper}). Hence, by Lemma \ref{lemma_of_poly_hyper_class_equivalent} we get the unique solution of (\ref{poly-hyperanalytic_RQ}) \(F\) in class \( \mathcal{R}_0^\mu(\overline{\partial}_q^{\,n}) \) and
\[ F(z)=\sum_{j=0}^{n-1} \overline{z}^j\Phi_j(z). \]

Case 3, \( \kappa<-1 \). Using Theorem \ref{hyper_theorem}, the necessary and sufficient condition for problem \( F_j^+(t)=G(t)F_j^-(t)+\check{f}_j(t) \) to be solvable in the class \( \mathcal{R}_0^\mu(\overline{\partial}_q) \) is
\[
\iint_{D^+} \frac{\overline{\partial}_q \check{f}^\omega_j(\zeta)}{X^+(\zeta)} W^{m-1}(\zeta) \, \mathrm{d}\xi \mathrm{d}\eta = 0, \quad m = 1, 2, \cdots, -\kappa - 1.
\]
Thus, the necessary and sufficient condition for Problem \ref{Jordan_curve_poly_hyper_RBVP} to be solvable in the class \( \mathcal{R}_0^\mu(\overline{\partial}_q^{\,n}) \) is
\[
\iint_{D^+} \frac{\overline{\partial}_q \check{f}^\omega_j(\zeta)}{X^+(\zeta)} W^{m-1}(\zeta) \, \mathrm{d}\xi \mathrm{d}\eta = 0, \quad m = 1, 2, \cdots, -\kappa - 1,\ j=0,1,\cdots,n-1.
\]
Furthermore, under the solvability condition, we can get that (\ref{poly-hyperanalytic_RQ}) has the unique solution
\[ F(z)=\sum\limits_{j=0}^{n-1} \overline{z}^j\Phi_j(z). \]
Thus, the theorem is proved.
\end{Proof}

\section{Riemann boundary value problem for meta-hyperanalytic functions}

In this section, we mainly solve the Riemann boundary value problem for meta-hyperanalytic functions. First, we give the definition of meta-hyperanalytic functions.

\begin{definition}
Let
\[ \overline{\partial}_{q,\lambda}^n:=( \overline{\partial}_q - \lambda )^{n} = \sum\limits_{k=0}^{n} \binom{n}{k} (-\lambda)^{n - k} \overline{\partial}_q^{\,k}, \] where \( \lambda\in\mathbb{D} \) is a hypercomplex number, \(n \) is a positive integer and \(\binom{n}{k}=\frac{n!}{(n-k)!k!}\). Let \( \Omega \) be an open subset of \( \mathbb{C} \). \( F\in C^n(\Omega) \) is called a meta-hyperanalytic function on \(\Omega\), if for each \(z\in \Omega\), \[ \overline{\partial}_{q,\lambda}^n F(z)=0.\] 
We will denote by \( \mathcal{M}_{q,\lambda}^n(\Omega)\) this class of functions.
\end{definition}

Let \( \Gamma \) be a closed Jordan curve. We call \(F\in \mathcal{M}_{q,\lambda}^n(\mathbb{C}\setminus\Gamma)\) a sectionally meta-hyperanalytic function with the jump \( \Gamma\), if the boundary values \( (\overline{\partial}_q^{\,k}F)^\pm(t) \) exist for all \( t\in\Gamma \) and \( k=0,1, \cdots,n-1 \).

The meta-hyperanalytic functions are related to the poly-hyperanalytic functions.

\begin{theorem}\label{meta_hyper_decomposition}
Let \(F\in\mathcal{M}_{q,\lambda}^n(\mathbb{C}\setminus\Gamma)\), then there exists unique \(E\in\mathcal{H}_q^n(\mathbb{C}\setminus\Gamma)\), such that
\[ F(z) = E(z)e^{\lambda \overline{z}},\quad z\in\mathbb{C}\setminus\Gamma,\]
where the exponential is understood in the sense of the Douglis algebra.
\end{theorem}
\begin{Proof}
Note that \(\partial e^{\lambda \overline{z}}=\frac{1}{2}(\lambda e^{\lambda \overline{z}}-\lambda e^{\lambda \overline{z}})=0 \) (since \(\partial(\overline{z}) = 0\)), then we have that
\begin{equation}\nonumber
\overline{\partial}_q^{\,n} \left( e^{-\lambda \overline{z}} F(z) \right) = \sum_{k=0}^n \binom{n}{k} \overline{\partial}_q^{\,k} e^{-\lambda \overline{z}} \cdot \overline{\partial}_q^{\,n-k} F(z)
= \sum_{k=0}^n \binom{n}{k} (-\lambda)^k e^{-\lambda \overline{z}} \cdot \overline{\partial}_q^{\,n-k} F(z)
= e^{-\lambda \overline{z}}\cdot \overline{\partial}_{q, \lambda}^n F(z).
\end{equation}
Hence \(E(z)=e^{-\lambda \overline{z}} F(z)\in\mathcal{H}_q^n(\mathbb{C}\setminus\Gamma)\), since \(\overline{\partial}_{q, \lambda}^n F(z)=0\). The uniqueness of \(E\) is obvious.
\end{Proof}

Using Theorem \ref{poly_hyper_decomposition}, for \(F \in \mathcal{M}_{q,\lambda}^n(\mathbb{C}\setminus\Gamma)\), there exist unique \(E_0,E_1,...,E_{n-1}\in\mathcal{H}_q(\mathbb{C}\setminus\Gamma)\), such that $$F(z)=\sum\limits_{k=0}^{n-1}\overline{z}^k E_k(z)e^{\lambda \overline{z}}.$$

The specific statement of the Riemann boundary value problem for meta-hyperanalytic functions is as follows.

\begin{problem}[Riemann boundary value problem for meta-hyperanalytic functions on a Jordan closed curve]\label{Jordan_curve_meta_hyper_RBVP}
Let \( \Gamma \) be a Jordan closed curve. Find a sectionally meta-hyperanalytic function \( F \) with the jump \( \Gamma \) such that
\begin{equation}\label{meta-hyperanalytic_RQ}
\begin{cases}
\begin{aligned}
F^+(t) &= G(t)F^-(t) + f_0(t), \\
(\overline{\partial}_q F)^+(t) &= G(t)(\overline{\partial}_q F)^-(t) + f_1(t), \\
\vdots & \\
(\overline{\partial}_q^{\,n-1} F)^+(t) &= G(t)(\overline{\partial}_q^{\,n-1} F)^-(t) + f_{n-1}(t),
\end{aligned}
\ t \in \Gamma,
\end{cases}
\end{equation}
where \( G,f_k \in H^\nu(\Gamma) \) are given for \( \nu\in (0,1] \), \( k=0,1,\cdots,n-1 \), \(G\) is a hypercomplex function whose complex part \(G_0\) never vanishes on \(\Gamma\).
\end{problem}

Next, to address this problem, we introduce some notations. For given \( \lambda \in \mathbb{D} \), let
\begin{equation}\nonumber
\mathbf{C}(\lambda) = \begin{pmatrix}
1 & 0 & 0 & \cdots & 0 \\
\lambda & 1 & 0 & \cdots & 0 \\
\lambda^2 & 2\lambda & 1 & \cdots & 0 \\
\vdots & \vdots & \vdots & \ddots & \vdots \\
\lambda^{n - 1} & \binom{n-1}{n-2}\lambda^{n - 2} & \binom{n-1}{n-3} \lambda^{n-3} & \cdots & 1
\end{pmatrix}.
\end{equation}
It is easy to see that \(\mathbf{C}(\lambda)\) is invertible and
\begin{equation}\nonumber
\mathbf{C}^{-1}(\lambda) = \begin{pmatrix}
1 & 0 & 0 & \cdots & 0 \\
-\lambda & 1 & 0 & \cdots & 0 \\
\lambda^2 & -2\lambda & 1 & \cdots & 0 \\
\vdots & \vdots & \vdots & \ddots & \vdots \\
(-\lambda)^{n - 1} & \binom{n-1}{n-2}(-\lambda)^{n - 2} & \binom{n-1}{n-3}(-\lambda)^{n - 3} & \cdots & 1
\end{pmatrix} = \mathbf{C}(-\lambda).
\end{equation}

For \( \mu\in(0,1] \), define the functions class
\begin{equation}\nonumber
\begin{split}
\mathcal{R}_m^\mu(\overline{\partial}_{q,\lambda}^n) =\{ F : &  F(z)=\sum\limits_{k=0}^{n-1}\overline{z}^kE_k(z)e^{\lambda\overline{z}} \text{ is a sectionally meta-hyperanalytic function with the jump }\Gamma, \\
&\mathrm{Ord}(E_k,\infty)\le m \text{ and } E_k \in H^\mu(\overline{D^+})\cap H^\mu(\overline{D^-}) \text{ for } k=0,1,\cdots,n-1\}.
\end{split}
\end{equation}

\begin{lemma}\label{lemma_of_meta_hyper_class_equivalent}
For \( F(z)=E(z)e^{\lambda\overline{z}}=\sum\limits_{k=0}^{n-1}\overline{z}^kE_k(z)e^{\lambda\overline{z}} \in \mathcal{M}_{q,\lambda}^n(\mathbb{C}\setminus\Gamma)\), \( F\in \mathcal{R}_0^\mu(\overline{\partial}_{q,\lambda}^n) \) and \( F \) satisfies the system of equations \((\ref{meta-hyperanalytic_RQ})\) if and only if \(E\in\mathcal{R}_0^\mu(\overline{\partial}_q^{\,n})\) and  \((\overline{\partial}_q^{\,k}E)^+(t)=G(t)(\overline{\partial}_q^{\,k}E)^-(t)+e^{-\lambda\overline{t} } \sum\limits_{j=0}^{k} \binom{k}{j} (-\lambda)^{k - j} f_{j}(t) \), for each \( k=0,1,\cdots,n-1 \).
\end{lemma}
\begin{Proof}
For \( F(z)=E(z)e^{\lambda\overline{z}}=\sum\limits_{k=0}^{n-1}\overline{z}^kE_k(z)e^{\lambda\overline{z}} \in \mathcal{M}_{q,\lambda}^n(\mathbb{C}\setminus\Gamma)\), applying \(\overline{\partial}_q^{\,k}\) \((k=0,1,\cdots,n-1)\) to both sides of the equation, we obtain that
\begin{equation}\label{n_equations_relate_meta_to_poly_hyper}
\begin{pmatrix}
F(z) \\
\overline{\partial}_q F(z) \\
\vdots \\
\overline{\partial}_q^{\,n - 1} F(z)
\end{pmatrix}
= e^{ \lambda \overline{z}} \mathbf{C}(\lambda)
\begin{pmatrix}
E(z) \\
\overline{\partial}_q E(z) \\
\vdots \\
\overline{\partial}_q^{\,n - 1} E(z)
\end{pmatrix}, \quad z\in \mathbb{C}\setminus\Gamma.
\end{equation}
From formula \((\ref{n_equations_relate_meta_to_poly_hyper})\) and the continuity of the exponential function in the Douglis algebra, it can be seen that \(F\) is a sectionally meta-hyperanalytic function with the jump \( \Gamma\) if and only if \(E\) is a sectionally poly-hyperanalytic function with the jump \( \Gamma\). Thus, we get that \( F\in \mathcal{R}_0^\mu(\overline{\partial}_{q,\lambda}^n) \)\(\Longleftrightarrow\)\( E \in \mathcal{R}_0^\mu(\overline{\partial}_q^{\,n}) \) by the definition of \(\mathcal{R}_0^\mu(\overline{\partial}_{q,\lambda}^n)\) and \(\mathcal{R}_0^\mu(\overline{\partial}_q^{\,n})\). Substitute (\ref{n_equations_relate_meta_to_poly_hyper}) into (\ref{meta-hyperanalytic_RQ}), we obtain that
\begin{equation}\nonumber
e^{ \lambda \overline{t} } \mathbf{C}(\lambda)
\begin{pmatrix}
E^+(t) \\
(\overline{\partial}_q E)^+(t) \\
\vdots \\
(\overline{\partial}_q^{\,n - 1} E)^+(t)
\end{pmatrix}
=
e^{ \lambda \overline{t}} \mathbf{C}(\lambda)
\begin{pmatrix}
G(t)E^-(t) \\
G(t)(\overline{\partial}_q E)^-(t) \\
\vdots \\
G(t)(\overline{\partial}_q^{\,n - 1} E)^-(t)
\end{pmatrix}
+\begin{pmatrix}f_0(t) \\ f_1(t) \\ \vdots \\f_{n-1}(t)\end{pmatrix}
, \quad t\in \Gamma.
\end{equation}
Let
\begin{equation}\label{hatf_hyper}
\hat{f}_k(t) = e^{-\lambda \overline{t} } \sum_{j=0}^{k} \binom{k}{j} (-\lambda)^{k - j} f_{j}(t), \quad k = 0, 1, \cdots, n-1,
\end{equation}
then we get the equivalent equation
\begin{equation}\label{equivalentmetaproblem_hyper}
\begin{cases}
\begin{aligned}
E^+(t) &= G(t)E^-(t) + \hat{f}_0(t), \\
(\overline{\partial}_q E)^+(t) &= G(t)(\overline{\partial}_q E)^-(t) + \hat{f}_1(t), \\
\vdots &\\
(\overline{\partial}_q^{\,n - 1} E)^+(t) &= G(t)(\overline{\partial}_q^{\,n - 1} E)^-(t) + \hat{f}_{n - 1}(t),
\end{aligned}
\end{cases}  t \in \Gamma,
\end{equation}
Then we have proved that \( F \) satisfies the system of equations \((\ref{meta-hyperanalytic_RQ})\) \(\Longleftrightarrow\)\( (\overline{\partial}_q^{\,k}E)^+(t)=G(t)(\overline{\partial}_q^{\,k}E)^-(t)+\hat{f}_k(t) \) for \( k=0,1,\cdots,n-1 \). Thus, the lemma holds.
\end{Proof}

Now, we present the solution to Problem \ref{Jordan_curve_meta_hyper_RBVP}. Let \( \kappa = \frac{1}{2\pi}\Delta_\Gamma \arg(G_0(t)) \), where \(G_0\) is the complex part of \(G\).

\begin{theorem}\label{solution_of_meta_hyper_RBVP}
Suppose that \(\Gamma\) is \(d\)-summable and \(d\in(1,2)\). If \( \nu > \dfrac{d}{2} \) and \( \mathrm{Dmh}\, \Gamma - 1 < \mu < \dfrac{2\nu - d}{2 - d} \) or \( \nu = 1 \) and \( \mathrm{Dmh}\, \Gamma - 1 < \mu < 1 \), then, when \( \kappa \ge 0 \), the general solution of Problem \textup{\ref{Jordan_curve_meta_hyper_RBVP}} in the class \( \mathcal{R}_0^\mu(\overline{\partial}_{q,\lambda}^n) \) is of the form \( F(z) =  e^{\lambda\overline{z}}\left(\sum\limits_{j=0}^{n-1}\overline{z}^j\Phi_j(z) + X(z)Q_{n-1,\kappa}(z,\overline{z})\right) \), where
\begin{equation}\nonumber
    \Phi_j(z) = h^\omega_j(z) + X(z) \iint_{D^+} e_z(\zeta) \left( \frac{\overline{\partial}_q h^\omega_j}{X} \right)(\zeta) \, \mathrm{d}\xi \mathrm{d}\eta,\quad j=0,1,\cdots,n-1,
\end{equation}
\(h_j\) is defined by \(h_j=\check{\hat{f}}_j\), \( X \) is defined as in \((\ref{hyper_X})\), and \( Q_{n-1,\kappa}(z,\overline{z}) \in \mathscr{P}_\kappa(W)+\overline{z}\mathscr{P}_\kappa(W)+\cdots+\overline{z}^{n-1}\mathscr{P}_\kappa(W)\). When \( \kappa = -1 \), \( F(z) = e^{\lambda\overline{z}}\sum\limits_{j=0}^{n-1}\overline{z}^j\Phi_j(z) \) is the unique solution of Problem \textup{\ref{Jordan_curve_meta_hyper_RBVP}} in the class \( \mathcal{R}_0^\mu(\overline{\partial}_{q,\lambda}^n) \). Finally, when \( \kappa < -1 \), the necessary and sufficient condition for Problem \textup{\ref{Jordan_curve_meta_hyper_RBVP}} to be solvable in the class \( \mathcal{R}_0^\mu(\overline{\partial}_{q,\lambda}^n) \) is
\begin{equation}\nonumber
    \iint_{D^+} \frac{\overline{\partial}_q h^\omega_j(\zeta)}{X^+(\zeta)} W^{m-1}(\zeta) \, \mathrm{d}\xi \mathrm{d}\eta = 0, \quad m = 1, 2, \cdots, -\kappa - 1,\ j=0,1,\cdots,n-1,
 \end{equation}
when this condition is satisfied, the function \( F(z) = e^{\lambda\overline{z}} \sum\limits_{j=0}^{n-1}\overline{z}^j\Phi_j(z) \) is also the unique solution of Problem \textup{\ref{Jordan_curve_meta_hyper_RBVP}} in the class \( \mathcal{R}_0^\mu(\overline{\partial}_{q,\lambda}^n)\).
\end{theorem}

\begin{Proof}
Consider the equations (\ref{equivalentmetaproblem_hyper}), where \( \hat{f}_j \) is defined by (\ref{hatf_hyper}), and, by Remark \ref{properties_of_Holder_class}, \( \hat{f}_j\in H^\nu(\Gamma) \) for \( j=0,1,\cdots,n-1 \). Using Theorem \ref{solution_of_poly_hyper_RBVP}, when \( \kappa \ge 0 \), the general solution of (\ref{equivalentmetaproblem_hyper}) in the class \( \mathcal{R}_0^\mu(\overline{\partial}_q^{\,n}) \) is of the form \( E(z) = \sum\limits_{j=0}^{n-1}\overline{z}^j\Phi_j(z) + X(z)Q_{n-1,\kappa}(z,\overline{z}) \), where
\begin{equation}\nonumber
    \Phi_j(z) = \left(\check{\hat{f}}_j\right)^\omega(z) + X(z) \iint_{D^+} e_z(\zeta) \left( \frac{\overline{\partial}_q \left(\check{\hat{f}}_j\right)^\omega}{X} \right)(\zeta) \, \mathrm{d}\xi \mathrm{d}\eta, \quad j=0,1,\cdots,n-1,
\end{equation}
\(h_k:=\check{\hat{f}}_k\) is defined by
\begin{equation}\nonumber
\begin{split}
h_k(t)=\check{\hat{f}}_k(t)&=\sum_{j=k}^{n-1}\dfrac{(-1)^{k+j}}{(j - k)!k!}\overline{t}^{j - k}\hat{f}_j(t)\\
&=\sum_{j=k}^{n-1}\dfrac{(-1)^{k+j}}{(j - k)!k!}\overline{t}^{j - k}\left(e^{- \lambda \overline{t}} \sum_{i=0}^{j} \binom{j}{i} (-\lambda)^{j - i} f_{i}(t)\right)\\
&=\sum_{j=k}^{n-1}\sum_{i=0}^{j}\dfrac{(-1)^{k+j}\binom{j}{i} (-\lambda)^{j - i}}{(j - k)!k!}e^{ - \lambda \overline{t}}\overline{t}^{j - k}f_{i}(t), \quad k=0,1,\cdots,n-1,
\end{split}
\end{equation}
\( X \) is defined as in (\ref{hyper_X}), and \( Q_{n-1,\kappa}(z,\overline{z}) \in \mathscr{P}_\kappa(W)+\overline{z}\mathscr{P}_\kappa(W)+\cdots+\overline{z}^{n-1}\mathscr{P}_\kappa(W)\). When \( \kappa = -1 \), \( E(z) = \sum\limits_{j=0}^{n-1}\overline{z}^j\Phi_j(z) \) is the unique solution of equations (\ref{equivalentmetaproblem_hyper}) in the class \( \mathcal{R}_0^\mu(\overline{\partial}_q^{\,n}) \). Finally, when \( \kappa < -1 \), the necessary and sufficient condition for equations (\ref{equivalentmetaproblem_hyper}) to be solvable in the class \( \mathcal{R}_0^\mu(\overline{\partial}_q^{\,n}) \) is
\begin{equation}\nonumber
    \iint_{D^+} \frac{\overline{\partial}_q h^\omega_j(\zeta)}{X^+(\zeta)} W^{m-1}(\zeta) \, \mathrm{d}\xi \mathrm{d}\eta = 0, \quad m = 1, 2, \cdots, -\kappa - 1,\ j=0,1,\cdots,n-1,
\end{equation}
when this condition is satisfied, the function \( E(z) = \sum\limits_{j=0}^{n-1}\overline{z}^j\Phi_j(z) \) is also the unique solution of equations (\ref{equivalentmetaproblem_hyper}) in the class \( \mathcal{R}_0^\mu(\overline{\partial}_q^{\,n})\). By Lemma \ref{lemma_of_meta_hyper_class_equivalent}, we obtain the conclusion in the theorem.
\end{Proof}

\begin{remark}
In the above theorem, the hypercomplex exponential \(e^{\lambda\overline{z}}\) is well-defined since the Douglis algebra admits a functional calculus for entire functions. The multiplication by \(e^{\lambda\overline{z}}\) preserves the H\"{o}lder continuity and the order at infinity of the components \(E_k\), thereby ensuring that the solution \(F\) belongs to the class \(\mathcal{R}_0^\mu(\overline{\partial}_{q,\lambda}^n)\).
\end{remark}

\section*{Competing interests}
No conflicts of interest are disclosed by the authors.

\section*{Data availability}
Data sharing not applicable to this article as no data sets were generated or analyzed during the current study.

\subsection*{ORCID}
\noindent
Juan Bory-Reyes: https://orcid.org/0000-0002-7004-1794\\
Fuli He: https://orcid.org/0000-0002-9395-545X

\end{document}